\documentclass[12pt]{article}
\usepackage{amsmath,amssymb,amsthm,cite}

\newtheorem{theorem}{Theorem}[section]
\newtheorem{thmy}{Theorem}

\newtheorem{corollary}[theorem]{Corollary}

\def\barr{\begin{array}}
\def\earr{\end{array}}

\title{A result on the Chermak-Delgado measure of a finite group}
\author{Marius T\u arn\u auceanu}
\date{August 15, 2024}

\begin{document}

\maketitle

\begin{abstract}
In this short note, we describe finite groups all of whose non-trivial cyclic subgroups have the same Chermak-Delgado measure.
\end{abstract}

{\small
\noindent
{\bf MSC2020\,:} Primary 20D30; Secondary 20D60, 20D99.

\noindent
{\bf Key words\,:} Chermak-Delgado measure, Chermak-Delgado lattice, subgroup lattice, TH-$p$-group.}

\section{Introduction}

Let $G$ be a finite group and $L(G)$ be the subgroup lattice of $G$. The \textit{Chermak-Delgado measure} of a subgroup $H$ of $G$ is defined by
\begin{equation}
m_G(H)=|H||C_G(H)|.\nonumber
\end{equation}Let
\begin{equation}
m^*(G)={\rm max}\{m_G(H)\mid H\leq G\} \mbox{ and } {\cal CD}(G)=\{H\leq G\mid m_G(H)=m^*(G)\}.\nonumber
\end{equation} Then the set ${\cal CD}(G)$ forms a modular, self-dual sublattice of $L(G)$, which is called the \textit{Chermak-Delgado lattice} of $G$. It was first introduced by Chermak and Delgado \cite{3}, and revisited by Isaacs \cite{5}. In the last years there has been a growing interest in understanding this lattice, especially for $p$-groups (see e.g. \cite{9}). We recall several important properties of the Chermak-Delgado measure:
\begin{itemize}
\item[$\bullet$] if $H\leq G$ then $m_G(H)\leq m_G(C_G(H))$, and if the measures are equal then $C_G(C_G(H))=H$;\newpage
\item[$\bullet$] if $H\in {\cal CD}(G)$ then $C_G(H)\in {\cal CD}(G)$ and $C_G(C_G(H))=H$;
\item[$\bullet$] the maximal member $M$ of ${\cal CD}(G)$ is characteristic and satisfies ${\cal CD}(M)={\cal CD}(G)$, while the minimal member $M(G)$ of ${\cal CD}(G)$ (called the \textit{Chermak-Delgado subgroup} of $G$) is characteristic, abelian and contains $Z(G)$.
\end{itemize}

In \cite{8}, the Chermak-Delgado measure of $G$ has been seen as a function
\begin{equation}
m_G:L(G)\longrightarrow\mathbb{N}^*,\, H\mapsto m_G(H),\, \forall\, H\in L(G),\nonumber
\end{equation}which has at least two distinct values if $G$ is non-trivial. We studied finite groups $G$ such that $m_G$ has exactly $k$ values, with an emphasis on the case $k=2$. Also, in \cite{4}, finite groups $G$ with $|{\cal CD}(G)|=|L(G)|-k$, $k=1,2$, have been determined. Note that a small $|{\rm Im}(m_G)|$ or a large ${\cal CD}(G)$ mean that many subgroups of $G$ have the same Chermak-Delgado measure. This constitutes the starting point of our discussion.

Our main result is stated as follows. By a \textit{TH-$p$-group} we will understand a $p$-group $G$ all of whose elements of order $p$ are central, that is $\Omega_1(G)\leq Z(G)$ (see e.g. \cite{1,2}).

\begin{theorem}
Let $G$ be a finite group and $C(G)^*$ be the set of non-trivial cyclic subgroups of $G$. If $m_G(H_1)=m_G(H_2)$ for all $H_1, H_2\in C(G)^*$, then $G$ is a TH-$p$-group with $\Omega_1(G)=Z(G)$. Moreover, if $|G|=p^n$, $\exp(G)=p^m$ and $|Z(G)|=p^k$, then $k\leq n-2m+2$.
\end{theorem}

Obviously, a finite abelian group as in Theorem 1.1 is an elementary abelian $p$-group. The smallest non-abelian examples are $Q_8$, $Q_8\times C_2$, $C_4\rtimes C_4$ for $p=2$ and $C_9\rtimes C_9$ for $p$ odd. More generally, we observe that all groups $Q_8\times C_2^n$ with $n\in\mathbb{N}$ and all groups $C_{p^2}\rtimes C_{p^2}$ with $p$ prime satisfy the hypothesis of Theorem 1.1.

Two particular cases of the above theorem are as follow.

\begin{corollary}
Let $G$ be a finite group. If any of the following two conditions holds
\begin{itemize}
\item[{\rm a)}] $m_G(H_1)=m_G(H_2)=m^*(G)$, for all $H_1, H_2\in C(G)^*$,
\item[{\rm b)}] $m_G(H_1)=m_G(H_2)$, for all non-trivial abelian subgroups $H_1, H_2$ of $G$,
\end{itemize}then either $G\cong C_p$ for some prime $p$ or $G\cong Q_8$.
\end{corollary}

Most of our notation is standard and will usually not be repeated here. Elementary notions and results on groups can be found in \cite{5,7}.
For subgroup lattice concepts we refer the reader to \cite{6}.\newpage

\section{Proof of the main results}

\bigskip\noindent{\bf Proof of Theorem 1.1.}
\medskip

Let $|G|=p_1^{n_1}\cdots p_r^{n_r}$ and $G_i\in{\rm Syl}_{p_i}(G)$, $i=1,...,r$. If $P_i$ is a cyclic subgroup of order $p_i$ which is contained in $Z(G_i)$, then $G_i\subseteq C_G(P_i)$ and so $m_G(P_i)$ is divisible by $p_i^{n_i+1}$. Since for $j\neq i$ the maximal power of $p_i$ in $m_G(P_j)$ is $p_i^{n_i}$, we infer that $m_G(P_i)\neq m_G(P_j)$. This shows that we must have $r=1$, i.e. $G$ is a $p$-group.

Next we observe that $Z(G)$ cannot contain elements of order $p^s$ with $s\geq 2$. Indeed, if $a$ is such an element, then $\langle a^p\rangle\neq 1$ and 
\begin{equation}
m_G(\langle a\rangle)>m_G(\langle a^p\rangle),\nonumber 
\end{equation}contradicting our hypothesis. Thus the common value of $m_G(H)$, $H\in C(G)^*$, is $p^{n+1}$.

Assume now that there is $H\leq G$ with $|H|=p$ and $H\nsubseteq Z(G)$. Then $C_G(H)\neq G$ and so
\begin{equation}
m_G(H)\leq p^n<p^{n+1},\nonumber
\end{equation}a contradiction. Consequently, $G$ is a TH-$p$-group with $\Omega_1(G)=Z(G)$.

Finally, let $b\in G$ with $o(b)=p^m$. Then both $\langle b\rangle$ and $Z(G)$ are contained in $C_G(\langle b\rangle)$, implying that $\langle b\rangle Z(G)\subseteq C_G(\langle b\rangle)$. This shows that
\begin{equation}
p^{m+k-1}=|\langle b\rangle Z(G)|\nonumber
\end{equation}divides $|C_G(\langle b\rangle)|$ and therefore $p^{2m+k-1}$ divides $p^{n+1}=m_G(\langle a\rangle)$. It follows that $2m+k-1\leq n+1$, i.e. 
\begin{equation}
k\leq n-2m+2,\nonumber
\end{equation}as desired.\qed

\bigskip\noindent{\bf Proof of Corollary 1.2.}
\medskip

Under the hypothesis of Theorem 1.1, we observe that if $k=1$, then $Z(G)$ is the unique subgroup of order $p$ of $G$. By (4.4) of \cite{7}, II, it follows that $G$ is either cyclic or a generalized quaternion $2$-group. Clearly, if $G$ is cyclic we must have $n=1$, that is $G\cong C_p$, where $p$ is a prime. If $G\cong Q_{2^n}$ for some integer $n\geq 3$, then $m=n-1$ and therefore the inequality $k\leq n-2m+2$ leads to $n\leq 3$. Thus $n=3$ and $G\cong Q_8$.

The proof is completed by the remark that each of conditions a) and b) implies $k=1$.\qed

\vspace*{3ex}\small

\hfill
\begin{minipage}[t]{5cm}
Marius T\u arn\u auceanu \\
Faculty of  Mathematics \\
``Al.I. Cuza'' University \\
Ia\c si, Romania \\
e-mail: {\tt tarnauc@uaic.ro}
\end{minipage}


\begin{thebibliography}{10}
\bibitem{1} D. Bubboloni, G. Corsi Tani, {\it $p$-groups with some regularity properties}, Ric. di Mat. {\bf 49} (2000), 327-339.
\bibitem{2} D. Bubboloni, G. Corsi Tani, {\it $p$-groups with all the elements of order $p$ in the center}, Algebra Colloq. {\bf 11} (2004), 181-190.
\bibitem{3} A. Chermak and A. Delgado, {\it A measuring argument for finite groups}, Proc. AMS {\bf 107} (1989), 907-914.
\bibitem{4} G. Fasol\u a and M. T\u arn\u auceanu, {\it Finite groups with large Chermak-Delgado lattices}, to appear in Bull. Aus. Math. Soc., DOI: https://doi.org/10.1017/S0004972722000806.
\bibitem{5} I.M. Isaacs, {\it Finite group theory}, Amer. Math. Soc., Providence, R.I., 2008.
\bibitem{6} R. Schmidt, {\it Subgroup lattices of groups}, de Gruyter Expositions in Ma\-the\-ma\-tics 14, de Gruyter, Berlin, 1994.
\bibitem{7} M. Suzuki, {\it Group theory}, I, II, Springer Verlag, Berlin, 1982, 1986.
\bibitem{8} M. T\u arn\u auceanu, {\it Finite groups with a certain number of values of the Chermak-Delgado measure}, J. Algebra Appl. {\bf 19} (2020), article ID 2050088.
\bibitem{9} A. Morresi Zuccari, V. Russo and C.M. Scoppola, {\it The Chermak-Delgado measure in finite $p$-groups}, J. Algebra {\bf 502} (2018), 262-276.
\end{thebibliography}
\end{document}